
\input amstex
\documentstyle{amsppt}
\magnification=1200
\catcode`\@=11
\redefine\logo@{}
\catcode`\@=13

\define \bn{\Bbb N}
\define \bz{\Bbb Z}
\define \bq{\Bbb Q}
\define \br{\Bbb R}
\define \bc{\Bbb C}

\define \M{{\Cal M}}

\define\La{{\Cal L}}

\define\geg{{\goth g}}
\define\rk{\text{rk}~}


\define\mult{\text{mult}}

\define\im{\text{im}}

\define\0o{{\overline 0}}
\define\1o{{\overline 1}}

\TagsOnRight

\document

\topmatter

\title
A Lecture about classification of Lorentzian Kac--Moody algebras of
the rank three
\endtitle

\author
Valeri A. Gritsenko and
Viacheslav V. Nikulin \footnote{Supported by
Grant of Russian Fund of Fundamental Research.
\hfill\hfill}
\endauthor

\address
University Lille 1,  
UFR de Mathematiques,  
F-59655 Villeneuve d'Ascq Cedex, 
France
\vskip1pt 
POMI, St. Petersburg, Russia 
\endaddress
\email
Valery.Gritsenko\@agat.univ-lille1.fr 
\endemail

\address
Deptm. of Pure Mathem. The University of Liverpool, Liverpool
L69 3BX, UK;
\vskip1pt 
Steklov Mathematical Institute,
ul. Gubkina 8, Moscow 117966, GSP-1, Russia
\endaddress
\email
vnikulin\@liv.ac.uk\ \ 
slava\@nikulin.mian.su
\endemail

\abstract
This is the announcement and partly the review of our results
about classification of Lorentzian Kac--Moody algebras of
the rank three. One of our results gives the classification
of Lorentzian Kac-Moody algebras with
denominator identity functions which are holomorphic 
reflective automorphic forms with respect to paramodular groups. 
Almost all of them were 
found in our paper ``Automorphic forms and Lorentzian Kac--Moody
algebras. Parts I, II'', Intern. J. Math. 9 (1998), 153--275. 
Our main result gives classification of all reflective 
meromorphic automorphic products with respect to paramodular groups.  

This is mainly a presentation of our lectures given during 
the Workshop: ``Applications to Quantum Field Theory'' ,
October 23 --- 27, 2000, in the Newton Mathematical 
Institute at Cambridge.
\endabstract

\rightheadtext
{Lorentzian Kac--Moody algebras}
\leftheadtext{V.A. Gritsenko and V.V. Nikulin}
\endtopmatter

\document

\head
0. Introduction
\endhead

In this paper we study Lorentzian Kac--Moody algebras with the
root lattice $S_t^\ast$, where
$S_t=U\oplus \langle 2t\rangle$, $t\in \bn$
(see definitions and notations below), and with the symmetry group
$$
\widehat{O}^+(L_t)=\{g\in O^+(L_t)\ |\ g \text{\ is trivial on
$L_t^\ast/L_t$}\}
$$
where $L_t=U\oplus S_t$. This case is especially interesting
because any hyperbolic lattice $S$ of the rank three and
representing zero has an 
equivariant embedding to $S_t$ or to its $m$-duals, 
$m|2t$. Any lattice $L$ of the signature
$(3,\,2)$ having a two dimensional isotropic sublattice has an
equivariant embedding to $L_t$ or to its $m$-duals, $m|2t$.
Moreover, one can consider only square-free $t$ here.
Thus, studying this case, we at the same time study 
Lorentzian Kac--Moody algebras with root lattices $S^\ast$
and symmetry groups $G\subset O^+(L)$ (of finite index)
where $S$ has an equivariant embedding to $S_t$ (or its $m$-duals), 
and $L$ has an equivariant embedding to $L_t$ (or its $m$-duals).

\head
1. Some general definitions on Lorentzian Kac--Moody algebras
\endhead

We consider the hyperbolic lattice $S_t:=U\oplus \langle 2t \rangle$,
$t\in \bn$. Here and in what follows a {\it lattice} means an integral
symmetric bilinear form which is usually denoted by $(\cdot\,,\,\cdot)$.
The lattice $U=\left( \matrix 0&-1\\-1&0\endmatrix\right)$, and
the lattice $\langle A\rangle$ is given by the symmetric
integral matrix $A$; we denote by $\oplus$ the orthogonal sum
of lattices. The lattice $U$ is hyperbolic, it has the signature
$(2,\,1)$. We also consider the lattice
$L_t=U\oplus S_t=U\oplus U\oplus \langle 2t \rangle$ of the
signature $(3,\,2)$.

In this situation a Lorentzian Kac--Moody algebra $\geg$ is given by
the {\it holomorphic automorphic form $\Phi (z)$ with respect to the
subgroup $G\subset O^+(L_t)$}
of finite index with the Fourier expansion of a very special form
$$
\Phi(z)=\sum_{w\in W}{\epsilon(w)
\Bigl(\exp{\left(-2\pi i (w(\rho),z)\right)}\ -\hskip-10pt
\bigr.\sum_{a\in S_t^\ast \cap \br_{++}\M}
{\bigl.m(a)\exp{\left(-2\pi i (w(\rho+a),z)\right)}\Bigr)}}
\tag{1.1}
$$
where {\it all coefficients $m(a)$ should be integral};
$W\subset O^+(S_t)$ is the {\it Weyl group};
$\epsilon:W\to \{\pm 1\}$ is its quadratic
character (e. g. $\epsilon=\det$ for usual cases);
$\M\subset \La (S_t)=V^+(S_t)/\br_{++}$ is
{\it a fundamental chamber of $W$;}
$\rho\in S_t\otimes \bq$ is the {\it Weyl vector} for the set
$P(\M)\subset (S_t)^\ast$ of {\it simple real roots.}
Additionally, the automorphic form
$\Phi(z)$ should be {\it reflective,} i. e. it should have zeros only in
{\it rational quadratic divisors} which are orthogonal to the set
$G(P(\M))$ of roots of $L_t$ where we always assume that $W\subset G$.

The automorphic form $\Phi$ automatically has the infinite product expansion
$$
\Phi(z)=\exp{(-2\pi i (\rho,z))}
\prod_{\alpha\in \Delta_+}{\Bigl(1-\exp{\left(-2\pi i (\alpha,z)\right)}
\Bigr)^{\mult(\alpha)}}.
\tag{1.2}
$$
where $\mult(\alpha)\in \bz$ and $\Delta_+\subset S_t^\ast $ is the
{\it set of positive roots of the algebra $\geg$.}

The sum part \thetag{1.1} defines the Lorentzian Kac--Moody
algebra $\geg$ by generators
$$
e_a,\,f_a,\,h_a,\ \ \
a\in P(\M)\bigcup
2P(\M)_{\1o}
\bigcup
\left( \bigcup_{a\in S_t^\ast \cap \br_{++}\M}{m(a)a}
\right)
\tag{1.3}
$$
(this is {\it the set of simple roots})
and defining relations (of Chevalley and Serre) in a usual way.
The algebra $\geg$ is graded by the dual lattice $S_t^\ast$:
$$
\geg(A)=\bigoplus_{\alpha\in S_t^\ast}{\geg_\alpha}=
\geg_0\bigoplus \left(\bigoplus_
{\alpha\in \Delta_+}\geg_\alpha\right) \bigoplus
\left(\bigoplus_{\alpha\in \Delta_+}
\geg_{-\alpha}\right),\ \ \ \geg_0=S_t\otimes \bc.
\tag{1.4}
$$
The product part \thetag{1.2} defines multiplicities
$$
\mult(\alpha):= \dim \geg_\alpha=\dim \geg_{\alpha,\0o}-
\dim \geg_{\alpha,\1o}.
$$  
In the standard terminology, the Lie algebra $\geg$ is
a {\it generalised Kac--Moody (or Borcherds) superalgebra} 
(an infinite-dimensional 
Lie superalgebra). See \cite{B1}, \cite{B2}, \cite{B3}, 
\cite{GN1}, \cite{GN5}, \cite{K}, \cite{N}, \cite{R} 
for details. Obviously, the product part \thetag{1.2}
can be also used to define the algebra $\geg$ since it defines
the sum \thetag{1.1}. The identity \thetag{1.1}=\thetag{1.2} is
called the {\it denominator identity} for the corresponding
algebra. We call the automorphic form \thetag{1.1}=\thetag{1.2}
the {\it denominator identity function.}

Below we explain more about the terminology above.

The Weyl group  $W\subset W(S_t)$ is a {\it reflection subgroup}
of the lattice $S_t$. It is generated by reflections in
some roots $\alpha \in S_t^\ast=Hom(S_t,\bz)$.
Here $\alpha \in S_t^\ast$ is called root if $\alpha^2>0$ and
$\alpha^2\,|\,2(S_t^\ast,\,\alpha)$. The same definition of roots
is used for any other lattice. The reflection group $W$ has
a fundamental chamber $\M$ in the hyperbolic space
$\La (S_t)=V^+(S_t)/\br_{++}$
where $V^+(S_t)$ is a half of the light cone
$V(S_t)=\{x\in S_t\otimes \br \ |\  x^2<0\}$. We denote by $P(\M)$
the set of orthogonal to $\M$ roots directed outwards of $\M$.
The set $P(\M)$ defines $\M$ and the group $W$ which
is generated by reflections in roots of $P(\M)$.

The Weyl vector $\rho \in S_t\otimes \bq$ is defined by the property.
$$
(\rho,\,\alpha)=-{\alpha^2\over 2}.
\tag{1.5}
$$
Existence of the Weyl vector is a very strong restriction on
$W$ and $P(\M)$.

The set $\Delta_+$ in \thetag{1.2} is called {\it the set of
positive roots}. It is union of the set of positive
imaginary roots $S\cap \overline{V^+(S)}$ and the set
$$
\{\alpha \in \bn\cdot W(P(\M))\ |\
\text {$\alpha$ is a root and\ \ } (\alpha,\ \M)<0\}
\tag{1.6}
$$
of positive real roots.

The variable $z$ in \thetag{1.1} and \thetag{1.2} belongs to
the complexified light cone 
$$
\Omega(V^+(S_t))=S_t\otimes \br+iV^+(S_t).
$$ 
It is canonically identified with the Hermitian symmetric domain
$$
\Omega(L_t)=\{\bc\omega \subset L_t\otimes \bc\ |\ (\omega,\,\omega)=0\ \and\
(\omega,\,\overline{\omega})<0 \}_0,
\tag{1.7}
$$
of the type IV as follows:
$z\in\Omega (V^+(S_t))$ defines the element
$\bc\omega_z\in \Omega(L_t)$
where $\omega_z =\big((z,z)/2\big)f_1+
f_{-1}\oplus z \in L_t\otimes \bc$
and $f_1$, $f_{-1}$ is the bases of the lattice $U$ with
the Gram matrix $U$ above. The coordinate $z$ is the 
affine coordinate of the connected symmetric domain 
$\Omega(L_t)$ in the neighbourhood of 
the cusp $f_1$ of the arithmetic group 
$O^+(L_t)$ which is the subgroup of $O(L_t)$ of index 
two acting on $\Omega(L_t)$. 

A holomorphic function $\Phi (z)$ on $\Omega (V^+(S))$ is an
{\it automorphic form of the weight $k$, $k\in \bz/2$,
with respect to a subgroup 
$G\subset O^+(L_t)$} of finite index if the function
$\widetilde{\Phi}(\lambda \omega_z)=\lambda^{-k}\Phi(z)$,
$\lambda \in \bc^\ast$, on the homogeneous cone 
$\widetilde{\Omega(L_t)}$ over $\Omega(L_t)$, 
satisfies the relation 
$\widetilde{\Phi}(g\omega)=\chi (g)\widetilde{\Phi}(\omega)$
for any $\omega \in \widetilde{\Omega(L_t)}_0$ and any
$g\in G$. Here $\chi$ is a finite character or a
multiplier system.

The rational quadratic divisor orthogonal to an
element $\alpha \in L_t\otimes \bq$ with $\alpha^2>0$ is equal to
$$
D_\alpha=\{\bc\omega \in \Omega (L_t)\ |\ (\omega,\,\alpha)=0\}.
\tag{1.8}
$$
Thus, we suppose that the divisor of $\Phi(z)$ is union of
rational quadratic divisors orthogonal to roots from
$G(P(\M))$.

In this definition, one can replace the lattice $S_t$ by
arbitrary hyperbolic (i.e. of signature $(m,1)$) lattice, and
the lattice $L_t$ by a lattice of signature $(m+1,2)$ with an 
isotropic primitive vector $e\in L$ such that $e^\perp = S$.
See \cite{B4}, \cite{GN5}.  
The lattices $S_t$ and $L_t$ above are especially
interesting (for the the rank three case) because they are maximal and
have maximal automorphism groups for square-free $t$.

The kernel $\widehat{O}^+(L_t)$ of the action of the
group $O^+(L_t)$ on the discriminant group
$A_{L_t}=(L_t^\ast/L_t)$ of the lattice $L_t$
is called the {\it extended paramodular group}. We want to describe
all automorphic forms of the form \thetag{1.1} with the properties
listed above which are automorphic with respect to the subgroup
$G\subset O^+(L_t)$ such that $G$ contains the extended
paramodular group $\widehat{O}^+(L_t)$.
It will give the classification of the
Lorentzian Kac--Moody algebras of the rank three
with the denominator identity function
\thetag{1.1}=\thetag{1.2} above which is an automorphic form with respect
to the extended paramodular group.

We have

\proclaim{Theorem 1.1} There are exactly 29 automorphic forms 
with respect to the extended paramodular group 
$\widehat{O}^+(L_t)$ defining Lorentzian Kac--Moody 
algebras of the rank three. They exist only for
$$
\split
t=&1\,(three),\ 2\,(seven),\ 3\,(seven),\ 4\,(seven),\ 8\,(one),\
9\,(one),\\
&12\,(one),\ 16\,(one),\ 36\,(one).
\endsplit
$$
where for each $t$ we show in brackets the number of forms.
All these 29 forms are given in Sect. 6 below.
\endproclaim

All these 29 forms were given in 
\cite{GN1}, \cite{GN2}, \cite{GN4},\cite{GN5} together with 
their sum and product expansions.  
Formally, three forms for $t=8,\,12,\,16$ are new but they 
coincide with some forms for $t=2$, $t=3$ and $t=4$ respectively 
after appropriate change of variable.  
Thus, in fact, they are not new. Perhaps, Theorem 1.1 is
the first result where a large class of Lorentzian Kac--Moody
algebras was classified. 
We describe all 29 forms of Theorem 1.1 in Sect. 6 below.

To construct automorphic forms of Theorem 1.1, we used the 
arithmetic lifting \cite{G1} of Jacobi forms which gives the Fourier 
expansion \thetag{1.1}. Also we used a variant of 
the Borcherds lifting \cite{B4} which we applied to Jacobi forms 
\cite{GN5}. It gives the infinite product 
expansion \thetag{1.2} of the forms. We consider this 
variant of the Borcherds lifting below.

\head
2. A variant of Borcherds products from \cite{GN5}
\endhead

We use a general result from \cite{GN5} 
which permits to construct many automorphic products of 
the form similar to \thetag{1.2} and which are automorphic with
respect to the extended paramodular group.

We use bases $f_2,\,f_{-2}$ for $U$ and $f_3$ for
$\langle 2t\rangle$. Together they give a bases $f_2,\,f_3,\,f_{-2}$
for the lattice $S_t=U\oplus \langle 2t \rangle$. The
dual lattice $S_t^\ast$ has the bases
$f_2,\,\widehat{f_3}=f_3/2t,\,f_{-2}$.
We denote
$$
\alpha=(n,l,m):=nf_2-l\widehat{f_3}+mf_{-2} \in S_t^\ast,
$$
where $\alpha^2=-2nm+{l^2\over 2t}$ and
$D(\alpha)=2t\alpha^2=-4tnm+l^2$ is the {\it discriminant or norm}  
of $\alpha$. 
For the dual lattice we usually use this form. Thus, we get the
lattice $S^\ast(2t)$.
We denote
$$
z=z_3f_2+z_2f_3+z_1f_{-2} \in \Omega(V^+(S_t)).
$$
Then
$$
\exp{\left(-2\pi i (\alpha,\,z)\right)}=q^nr^ls^m
$$
where $q=exp{(2\pi i z_1)}$, $r=\exp{(2\pi i z_2)}$,
$s=\exp{(2\pi i z_3)}$.

In \cite{GN5} a variant of the Borcherds lifting was proved. We formulate
that in Theorem 2.1 below.

Let
$$
\phi_{0,t}(\tau,z)=\sum_{k,l\in \bz} f(k,l)q^kr^l\in
J_{0,t}^{nh}
\tag{2.1}
$$
be a nearly holomorphic Jacobi form of weight $0$
and index $t$ (i.e. $n$ might be negative in the Fourier expansion),
where
$q=\exp{(2\pi i \tau)}$, $\im \tau>0$, $r=\exp{(2\pi i z)}$.
It is automorphic with respect to the Jacobi group
$H(\bz) \rtimes SL_2(\bz)$ where
$H(\bz)$ is the integral Heisenberg group which is the central extension
$$
0 \to \bz \to H(\bz)\to \bz\times \bz\to 0.
$$
Let
$$
\phi^{(0)}_{0,t}(z)=\sum_{l\in \bz} f(0,l)\,r^l
\tag{2.2}
$$
be the $q^0$-part of $\phi_{0,t}(\tau,\,z)$.
The Fourier coefficient $f(k,l)$ of $\phi_{0,t}$ depends only on
the norm  $4tk-l^2$ of $(k,l)$ and $l$ mod $2t$.
Moreover, $f(k,l)=f(k,-l)$.
From the definition of nearly holomorphic forms, it follows
that the norm $4tk-l^2$ of indices  of non-zero Fourier coefficients
$f(k,l)$ are bounded from bellow.

\proclaim{Theorem 2.1 (\cite{GN5, Part II, Theorem 2.1})} 
Assume that the Fourier coefficients $f(k,l)$ of a 
Jacobi form $\phi_{0,t}$ 
from \thetag{2.1} are integral.
Then the infinite product
$$
B_\phi(z)=
q^{A}r^Bs^{C}
\prod\Sb n,l,m\in \bz\\
\vspace{0.5\jot}
(n,l,m)>0\endSb
 (1-q^nr^ls^m)^{f(nm,l)},
\tag{2.3}
$$
where
$$
A=\frac{1}{24}\sum_{l}f(0,l),\quad
B=\frac{1}{2}\sum_{l>0}lf(0,l),\quad
C=\frac{1}{4t}\sum_{l}l^2f(0,l),
$$
and $(n,l,m)>0$ means that either $m>0$ or $m=0$ and $n>0$ or
$m=n=0$ and $l<0$, defines  a meromorphic modular form of weight
$\frac{f(0,0)}2$ with respect to $\widehat{O}^+(L_t)$ with a character
(or a multiplier system if the weight is half-integral)
induced by $v_\eta^{24A}\times v_H^{4tC}$.   
All divisors of $B_\phi(z)$
are the rational quadratic divisors orthogonal to
$\alpha=(a,b,c)>0$ of the discriminant $D=-4tac+b^2$
with multiplicities
$$
m_{ac,b}=\sum_{n>0}f(n^2ac,nb).
\tag{2.4}
$$
Moreover
$$
B_\phi(V_t(z))=(-1)^{D} B_\phi(z)
\qquad\text{with}\quad
D=\sum\Sb k<0 \\
\vspace{0.5\jot} l \in \bz \endSb
\sigma_1(-k)f(k,l)
\tag{2.5}
$$
where
$\sigma_1(k)=\sum_{d|k}d$. See details in \cite{GN5, Part II, Theorem 2.1}.
\endproclaim

All 29 automorphic forms of Theorem 1.1 are given by some automorphic
products of Theorem 2.1. Thus, to give all these 29 automorphic forms,
we should give the corresponding 29 Jacobi forms. We give all these
forms in Sect. 6, Table 2.  

More generally, we describe all reflective automorphic forms
which are given by the Theorem 2.1 where a meromorphic automorphic form
on the domain $\Omega(L_t)$ is called {\it reflective} if its
divisor is union of rational quadratic divisors orthogonal to roots of
$L_t$. We remind that an element $\alpha \in L_t$ is called
{\it root} if $\alpha^2\,|\,2(\alpha,\,L_t)$ and $\alpha^2>0$.
One can prove that the infinite product $B_\phi$ given by 
a Jacobi form $\phi=\phi_{0,t}$ is reflective if and only if 
each non-zero Fourier coefficient $f(k,l)$ of $\phi_{0,t}$ 
with negative norm $4tk-l^2<0$ satisfies 
$$
4tk-l^2\,|\,(4t,\,2l).
$$
Let us {\it denote by $RJ_t$} the space of all Jacobi forms of the 
index $t$ of Theorem 2.1 which give reflective automorphic products. 
It is natural to call these Jacobi forms {\it reflective} either.
The space $RJ_t$ of all reflective Jacobi forms
is a free $\bz$-module with respect to addition.

\proclaim{Main Theorem 2.2} For $t\in \bn$ the space $RJ_t$ of
reflective Jacobi forms of the index $t$ is not trivial 
(i.e. is not equal to zero) if and only if $t$ is equal to
$$
\split
&1(2),\,2(3),\,3(3),\,4(3),\,5(3),\,6(4),\,7(2),\,8(3),\,9(3),\,
10(3),\,11(1),\,12(4),\,13(2),\\
&14(3),\,15(2),\,16(2),\,17(1),\,18(3),\,20(3),\,21(3),
\,22(1),\,24(2),\,25(1),\,26(1),\\
&28(1),\,30(3),\,33(1),\,34(2),\,36(3),\,39(2),\,42(1),\,45(1),\,
48(1),\,63(1),\,66(1).
\endsplit
$$
where in brackets we also show the rank of the corresponding
$\bz$-module $RJ_t$ of reflective Jacobi forms.

In Sect. 5, Table 1 below we give the bases of the module $RJ_t$ for 
$t=1$, $2$, $3$, $4$, $8$, $9$, $12$, $16$, $36$ 
when the subspace $RJ_t$ also contains a Jacobi form which gives
the denominator identity for a Lorentzian Kac--Moody algebra
(i.e. it gives the forms of Theorem 1.1). For all other $t$
the list of all Jacobi forms of the theorem
will be given in the forthcoming publication.
\endproclaim

All automorphic forms of Theorem 1.1 are characterised by the
property that they have multiplicity 1 for
components of their divisors. So, it is not hard to find all
these forms from the corresponding full Table of reflective
Jacobi forms of Main Theorem 2.2 since the theorem 2.1 also
gives the multiplicities of divisors. See Table 1 for
$t=1$, $2$, $3$, $4$, $8$, $9$, $12$, $16$, $36$.

Potentially, Main Theorem 2.2 contains information about all
reflective automorphic forms with infinite product expansion
of the type of Theorem 2.1 for all equivariant sublattices
$L\subset L_t$ of finite index.
Here {\it equivariant} means
that $O(L)\subset O(L_t)$; in particular, every root of $L$ is
multiple to a root of $L_t$. If $L$ has a reflective automorphic
form $\Phi$ with respect to $O(L)$ with an infinite product 
expansion, then its symmetrization
$$
\prod_{g\in O(L)\backslash O(L_t)}{g^\ast\Phi}
$$
is a reflective automorphic form with an infinite product expansion
for the lattice $L_t$.

Thus, potentially, Main Theorem 2.2 contains
important information about reflective automorphic forms with
infinite products and about automorphic forms of denominator identities
of Lorentzian Kac--Moody algebras with the lattices $L$ instead of $L_t$
and the corresponding hyperbolic lattices $S=S_t\cap L$ instead of $S_t$.
Moreover, one can possibly consider more general class of
Lie algebras for which reflective forms of Main Theorem 2.2 may
give some denominator identities. Results from \cite{GaL} and 
\cite{KW} give a hope.

\head
3. The proof of Main Theorem 2.2 and
reflective hyperbolic lattices.
\endhead

To classify finite-dimensional semi-simple or 
affine Lie algebras, one needs to 
classify corresponding finite or affine root systems.  
To prove Main Theorem 2.2 one needs description of appropriate 
hyperbolic root systems.  

Let $S$ be a hyperbolic (i.e. of the signature $(m,1)$) lattice,
$W(S)$ its reflection group and $\M\subset V^+(S)/\br_{++}$ its
fundamental chamber and $A(\M)$ is the symmetry group of the fundamental
chamber. Thus, $O^+(S)=W(S)\rtimes A(\M)$. The lattice $S$ is called
{\it reflective} if $\M$ has a {\it generalised Weyl vector $\rho\in
S\otimes \bq$}. It means that $\rho\not=0$ and the orbit $A(\M)(\rho)$
is finite. A reflective lattice is called
{\it elliptically reflective} if it has a generalised
Weyl vector $\rho$ with $\rho^2<0$. It is called {\it parabolically
reflective} if it is not elliptically reflective but has
a generalised Well vector $\rho$ with $\rho^2=0$. It is called
{\it hyperbolically reflective} if it is not elliptically or
parabolically reflective, but it has a generalised Weyl vector
$\rho$ with $\rho^2>0$.

Suppose that $B_\phi(z)$ is a reflective automorphic form
of Theorem 2.1.
The inequality $(n,m,l)>0$ of Theorem 2.1 is a variant of
choosing a fundamental chamber $\M$ of $W(S_t)$. It follows that
if the vector $\rho=(A,B,C)$ is not zero, then it defines a generalised
Weyl vector for $A(\M)$. The vector $\rho$ is invariant with
respect to the group $\widehat{A}(\M)=A(\M)\cap \widehat{O}(S_t)$ which
has finite index in $A(\M)$.

If the form $B_\phi(z)$ has a zero Weyl vector $\rho=(A,B,C)$, one 
can change it by other form which will have a non-zero Weyl vector,
considering reflections in roots. Thus, we get 

\proclaim{Lemma 3.1} If the space $RJ_t$ of reflective Jacobi forms
is not zero, then the lattice $S_t$ is reflective.
\endproclaim

It is interesting that the space $RJ_t$ may really have a
Jacobi form with zero Weyl vector $\rho=(A,B,C)$. It happens for
$t=6$ and $t=12$ when $\rk RJ_t=4$.

In \cite{N5}, for the rank three, 
all maximal reflective hyperbolic lattices were classified. 
More generally, all reflective hyperbolic lattices $S$ with
the square-free determinant $\det (S)$ were classified for the 
rank three. In the 
same paper we gave estimate for invariants of any reflective
hyperbolic lattice of the rank three. Using this estimate, one
can obtain classification of reflective hyperbolic lattices 
$S_t=U\oplus \langle2t \rangle$ for all $t\in\bn$ (see also 
calculations in \cite{N4}). As a result, we have

\proclaim{Theorem 3.2} The lattice $S_t=U\oplus \langle 2t \rangle$
is reflective for the
following and the only following $t\in \bn$,
where in brackets we put the
type of reflectivity of the lattice, (e) for elliptic, (p) for parabolic
and (h) for the hyperbolic type:
$$
\split
t\,=\,&1\ -\ 22\,(e),\ 23\,(h),\ 24\ -\ 26\,(e),\ 28\,(e),\
29\,(h),\ 30\,(e),\ 31\,(h),\\
&33\,(e),\ 34\,(e),\ 35\,(h),\ 36\,(e),\ 37\,(h),\ 38\,(h),\ 39\,(e),\
40\,(h),\ 42\,(e),\\
&44\,(h),\ 45\,(e),\ 46\,(h),\ 48\,(h),\ 49\,(e),\ 50\,(e),\ 52\,(e),\
55\,(e),\ 56\,(h),\\
&57\,(h),\ 60\,(h),\ 63\,(h),\ 66\,(e),\ 70\,(h),\ 72\,(h),\
78\,(h),\ 84(h),\ 90\,(h),\\ 
&100\,(h),\ 105\,(h).
\endsplit
$$
\endproclaim

To prove Main Theorem 2.2 and to find the bases of $RJ_t$,
one needs to analyse only the $t$ of Theorem 3.2.
These can be done using the generators of the graded ring of 
the weak Jacobi forms with integral Fourier coefficients found 
in \cite{G2} and \cite{G3}. We give these bases 
for $t=1$, $2$, $3$, $4$, $8$, $9$, $12$, $16$, $36$ in
Sect. 5, Table 1 below.

One can also use arguments
which we give below for the proof of Theorem 1.1.

\smallpagebreak

We mention that in \cite{N5} the classification of all
reflective hyperbolic lattices with square-free determinant 
(and their duals) is given for the rank three.  
In particular, it contains classification
of all maximal reflective hyperbolic lattices for the rank three. 
E.g. this classification contains $122\,(e) + 66\,(h)$ main 
hyperbolic lattices with square-free determinant, only 
$23\,(e)+11\,(h)$ of them represent $0$ and are given in 
Theorem 3.2. Moreover, in \cite{N5}, there are estimates on
invariants of all reflective hyperbolic lattices for the rank three. 
Thus, potentially, one can use these results for classification
of all reflective automorphic forms with infinite products of
the type (see below) of Main Theorem 2.2, for the rank three.  

\head
4. Proof of Theorem 1.1 and reflective hyperbolic lattices
with Weyl vector.
\endhead

Here we sketch the proof of Theorem 1.1 to emphasise importance of
reflective hyperbolic lattices with Weyl vector.

For this case, the fundamental polyhedron $\M$ and the set $P(\M)$
of roots orthogonal to $\M$ have the Weyl vector $\rho$ (satisfying
\thetag{1.5}). They are all
invariant with respect to the group $\widehat{A}(\M)$
(we use notation $\widehat{G}=G\cap \widehat{O}^+(L_t)$).
For all reflective lattices $S_t$ the
fundamental polyhedron $\M_0$ for the full reflection group $W(S_t)$ of
the lattice $S_t$ can be calculated and is known (for
$t=1$, $2$, $3$, $4$, $8$, $9$, $12$, $16$, $36$
these calculations are presented in Table 1).
Thus, the fundamental
chamber $\M$ is composed from the known polyhedron $\M_0$
by some reflections. Using this information, we can find all possible
$\M$, $P(\M)$, $\rho$ and predict the divisor
of the reflective automorphic form $\Phi(z)$. Looking at the list
of 29 forms of Theorem 1.1, one can see that one of the forms (for
the corresponding $t$) of Theorem 1.1 has the same divisor. By
Koecher principle, the form $\Phi(z)$ is equal to that form.

 Similar arguments can be used to classify all {\it reflective
meromorphic automorphic forms with infinite product}
of the form similar to \thetag{1.2} and with a generalised Weyl vector
$\rho$.  Like the product \thetag{1.2},
this product is related with a reflection subgroup $W\subset W(S_t)$,
its fundamental chamber $\M$, the set $P(\M)$ of orthogonal roots
to $\M$ (they define $\Delta_+$) and a generalised Weyl vector
$\rho\in S_t\otimes \bq$, i. e.  $\rho\not=0$, the orbit
$A(\M)(\rho)$ is finite and
$W\rtimes A(\M)$ has finite index in $O(S_t)$. The function
$\mult(\alpha)$, $\alpha \in \Delta_+$, should be integral and
invariant with respect to $\widehat{A}(\M)$.
The product should converge in a neighbourhood of the cusp $im (z)^2<<0$.
All definitions are the same. We have

\proclaim{Theorem 4.1} All reflective meromorphic automorphic forms
with respect to the extended paramodular group and
with infinite product, where $\rho$ is
a non-zero generalised Weyl vector, belong to the list of
Main Theorem 2.2.
\endproclaim

Applying to the forms of Theorem 4.1 reflections with respect to roots
in $S_t$, one can get some automorphic forms with zero Weyl vector and with
infinite product. They appear only for $t=6$ or $t=12$.

\head
5. The list of all reflective Jacobi forms from $RJ_t$ for
$t=1$, $2$, $3$, $4$, $8$, $9$, $12$, $16$, $36$.
\endhead

For these $t$ we give the bases
$\xi_{0,t}^{(1)},\,\dots\,\xi_{0,t}^{(rk)}$ of the
$\bz$-module $RJ_t$ showing the
leading parts of their Fourier expansions which define the
Jacobi form uniquely. We give all their Fourier coefficients
with the negative norm (up to equivalence); the corresponding
negative norm is shown in brackets $[\cdot]$. We also give
expressions of $\xi_{0,t}^{(i)}$ 
using basic Jacobi forms. In these formulae 
$E_4=E_4(\tau)$ and $\Delta_{12}=\Delta(\tau)$ are the Eisenstein series
of weight $4$
and the Ramanujan function of weight $12$ for $SL_2(\bz )$,
$E_{4,m}$ ($m=1,\,2,\,3$) are  Eisenstein-Jacobi
series of weight 4 and index $m$ (see \cite{EZ}), and 
 $\phi_{0,1}$, $\phi_{0,2}$, $\phi_{0,3}$, $\phi_{0,4}$
are the four generators of the graded ring of the weak Jacobi forms
of weight zero with integral Fourier coefficients (see \cite{G2} and 
\cite{GN5}).

We give the set $\overline{R}$ of primitive roots in
$S_t^\ast$ up to equivalence (up to the action of the group
$\pm \widehat{O}(S_t)$). Up to this equivalence, a root
$\alpha=(n,\,l,\,m)$ is
defined by its norm $-2t\alpha^2=-4nm+l^2$ and $l\mod 2t$.
We also give the matrix
$$
Mul(\overline{R},\,\xi)=mul(\gamma_i,\,\xi_{0,t}^{(j)}),
$$
where $mul(\gamma_i,\,\xi_{0,t}^{(j)})$ is the multiplicity of
the form $\xi_{0,t}^{(j)}$ in
a rational quadratic divisor which is orthogonal to the root from
the equivalence class $\gamma_i\in \overline{R}$.

We give the set $P(\M_0)$ of primitive roots in $S_t^\ast$ which is
orthogonal to the fundamental chamber $\M_0$ of the
reflection group $W(S_t)$ (this is equivalent to
the ordering $(n,l,m)>0$ used in Theorem 2.1), and their Gram matrix
$$
G(P(\M_0))=2t((\alpha,\,\beta)),\ \alpha,\,\beta \in P(\M_0).
$$
Thus, we identify the dual lattice $S^\ast$ with the integral lattice
$S^\ast(2t)=U(2t)\oplus \langle 1 \rangle$ to make it integral.

All these data are given in Table 1 below.

\vskip40pt

\centerline{\bf Table 1. The spaces $RJ_t$ for
$t=1$, $2$, $3$, $4$, $8$, $9$, $12$, $16$, $36$.}

\vskip20pt

\centerline{Case $t=1$.}

\vskip10pt

The space $RJ_1$ has the bases
$$
\xi_{0,1}^{(1)}=\phi_{0,1}=
$$
$$
=(r[-1] + 10 + r^{-1}[-1]) + O(q); 
$$
$$
\xi_{0,1}^{(2)}={E_4}^2{E_{4,1}}/{\Delta_{12}}-57\phi_{0,1}=
$$
$$
= q^{-1}[-4] +
(r^2[-4]-r[-1]+60-r^{-1}[-1]+r^{-2}[-4]) + O(q).
$$
We have $\overline{R}=\overline{P(\M_0)}$ and
$$
P(\M_0)=
\pmatrix
{1}&{2}&{0}\cr
{0}&{-1}&{0}\cr
{-1}&{0}&{1}\cr
\endpmatrix \equiv
\left[\matrix
4,\,\overline{0}\cr
1,\,\overline{1}\cr
4,\,\overline{0}\cr
\endmatrix\right];
\hskip10pt
Mul(P(\M_0),\,\xi)=
\pmatrix
0&1\cr
1&0\cr
0&1\cr
\endpmatrix;
$$
$$
G(P(\M_0))=
\pmatrix
{4}&{-2}&{-2}\cr
{-2}&{1}&{0}\cr
{-2}&{0}&{4}\cr
\endpmatrix.
$$

\vskip20pt

\vskip20pt

\centerline{Case $t=2$.}

\vskip10pt

The space $RJ_2$ has the bases
$$
\xi_{0,2}^{(1)}={\phi_{0,2}}=
$$
$$
=(r[-1] + 4 + r^{-1}[-1]) + O(q);
$$
$$
\xi_{0,2}^{(2)}=({\phi_{0,1}})^2-21{\phi_{0,2}}=
$$
$$
=(r^2[-4]- r[-1]+18-r^{-1}[-1]+r^{-2}[-4]) + O(q);
$$
$$
\xi_{0,2}^{(3)}={E_4}^2{E_{4,2}}/{\Delta_{12}}-
14({\phi_{0,1}})^2+216{\phi_{0,2}}=
$$
$$
= q^{-1}[-8] + 24 + O(q).
$$
We have $\overline{R}=\overline{P(\M_0)}$ where
$$
P(\M_0)=
\pmatrix
{1}&{2}&{0}\cr
{0}&{-1}&{0}\cr
{-1}&{0}&{1}\cr
\endpmatrix
\equiv
\left[\matrix
4,\,\overline{2}\\
1,\,\overline{1}\\
8,\,\overline{0}\\
\endmatrix\right];
\hskip20pt
Mul(P(\M_0),\,\xi)=
\pmatrix
0&1&0\cr
1&0&0\cr
0&0&1\cr
\endpmatrix;
$$
$$
G(P(\M_0))=
\pmatrix
{4}&{-2}&{-4}\cr
{-2}&{1}&{0}\cr
{-4}&{0}&{8}\cr
\endpmatrix.
$$

\vskip20pt

\centerline{Case $t=3$.}

\vskip10pt

The space $RJ_3$ has the bases
$$
\xi_{0,3}^{(1)}={\phi_{0,3}}=
$$
$$
=(r[-1] + 2 + r^{-1}[-1]) + O(q);
$$
$$
\xi_{0,3}^{(2)}={\phi_{0,1}}{\phi_{0,2}}-15{\phi_{0,3}}=
$$
$$
=(r^2[-4] -r[-1] + 12 -r^{-1}[-1]+ r^{-2}[-4]) + O(q);
$$
$$
\xi_{0,3}^{(3)}={E_4}^2{E_{4,3}}/{\Delta_{12}}-2({\phi_{0,1}})^3+
33{\phi_{0,1}}{\phi_{0,2}}+90{\phi_{0,3}}=
$$
$$
= q^{-1}[-12] + 24 + O(q).
$$
We have $\overline{R}=\overline{P(\M_0)}$ where
$$
P(\M_0)=
\pmatrix
{1}&{2}&{0}\cr
{0}&{-1}&{0}\cr
{-1}&{0}&{1}\cr
\endpmatrix
\equiv
\left[\matrix
4,\,\overline{2}\\
1,\,\overline{1}\\
12,\,\overline{0}\\
\endmatrix\right];
\hskip20pt
Mul(P(\M_0),\,\xi)=
\pmatrix
0&1&0\cr
1&0&0\cr
0&0&1\cr
\endpmatrix;
$$
$$
G(P(\M_0))=
\pmatrix
{4}&{-2}&{-6}\cr
{-2}&{1}&{0}\cr
{-6}&{0}&{12}\cr
\endpmatrix.
$$

\vskip20pt

\centerline{Case $t=4$.}

\vskip10pt

The space $RJ_4$ has the bases
$$
\xi_{0,4}^{(1)}={\phi_{0,4}}=
$$
$$
=(r[-1] + 1 + r^{-1}[-1]) + O(q);
$$
$$
\xi_{0,4}^{(2)}=({\phi_{0,2}})^2-9{\phi_{0,4}}=
$$
$$
=(r^2[-4] - r[-1] + 9 -r^{-1}[-1]+ r^{-2}[-4]) + O(q);
$$
$$
\xi_{0,4}^{(3)}={E_4}{E_{4,1}}{E_{4,3}}/{\Delta_{12}}-
2({\phi_{0,1}})^2{\phi_{0,2}}+20{\phi_{0,1}}{\phi_{0,3}}+16{\phi_{0,4}}=
$$
$$
= q^{-1}[-16] + 24 + O(q).
$$
We have $\overline{R}=\overline{P(\M_0)}$ where
$$P(\M_0)=
\pmatrix
{1}&{2}&{0}\cr
{0}&{-1}&{0}\cr
{-1}&{0}&{1}\cr
\endpmatrix
\equiv
\left[\matrix
4,\,\overline{2}\\
1,\,\overline{1}\\
16,\,\overline{0}\\
\endmatrix\right];
\hskip20pt
Mul(P(\M_0),\,\xi)=
\pmatrix
0&1&0\cr
1&0&0\cr
0&0&1\cr
\endpmatrix;
$$
$$
G(P(\M_0))=
\pmatrix
{4}&{-2}&{-8}\cr
{-2}&{1}&{0}\cr
{-8}&{0}&{16}\cr
\endpmatrix.
$$

\vskip20pt

\centerline{Case $t=8$.}

\vskip10pt

The space $RJ_8$ has the bases
$$
\xi_{0,8}^{(1)}=({\phi_{0,2}})^2{\phi_{0,4}}-
{\phi_{0,2}}({\phi_{0,3}})^2-({\phi_{0,4}})^2=
$$
$$
\split
=&(2r[-1] - 1 + 2r^{-1}[-1]) +(-r^6[-4] - 2r^5 + 4r^4 -
4r^3 + r^2 + 6r - 8+\\
&6r^{-1}+ r^{-2} - 4r^{-3} + 4r^{-4} - 2r^{-5} - r^{-6}[-4])q + O(q^2);
\endsplit
$$
$$
\xi_{0,8}^{(2)}=\phi_{0,2}(\tau,2z)=
{\phi_{0,1}}{\phi_{0,3}}{\phi_{0,4}}+{\phi_{0,2}}({\phi_{0,3}})^2-
2({\phi_{0,2}})^2{\phi_{0,4}}-2({\phi_{0,4}})^2=
$$
$$
=(r^2[-4]+4+r^{-2}[-4])+
(r^6[-4]-8r^4-r^2+16
-r^{-2}-8r^{-4}+r^{-6}[-4])q+O(q^2);
$$
$$
\split 
\xi_{0,8}^{(3)}=&{E_4}{E_{4,3}}\Big({E_{4,2}}{\phi_{0,3}}-
{E_{4,1}}{\phi_{0,4}}\Big)/{\Delta_{12}}-
3({\phi_{0,1}})^2({\phi_{0,3}})^2+\\
&2({\phi_{0,1}})^2{\phi_{0,2}}{\phi_{0,4}}+
8{\phi_{0,1}}{\phi_{0,3}}{\phi_{0,4}}-16({\phi_{0,4}})^2=
\endsplit 
$$
$$
\split
=&q^{-1}[-32] + 24 +\\
&(8r^6[-4] + 256r^5 + 2268r^4 +9472r^3+23608r^2 + 39424r + 46812 +\\
&39424r^{-1} + 23608r^{-2} + 9472r^{-3} +
2268r^{-4} + 256r^{-5} + 8r^{-6}[-4])q + O(q^2).
\endsplit
$$
We have $\overline{R}=
\overline{P(\M_0)}$ where
$$
P(\M_0)=\pmatrix
{1}&{2}&{0}\cr
{0}&{-1}&{0}\cr
{-1}&{0}&{1}\cr
{1}&{6}&{1}\cr
\endpmatrix
\equiv
\left[\matrix
4,\,\overline{2}\\
1,\,\overline{1}\\
32,\,\overline{0}\\
4,\,\overline{6}
\endmatrix\right];
\hskip20pt
Mul(P(\M_0),\,\xi)=
\pmatrix
0  &1  &0\cr
2  &1  &0\cr
0  &0  &1\cr
-1 &1  &8
\endpmatrix;
$$
$$
G(P(\M_0))=
\pmatrix
{4}&{-2}&{-16}&{-4}\cr
{-2}&{1}&{0}&{-6}\cr
{-16}&{0}&{32}&{0}\cr
{-4}&{-6}&{0}&{4}\cr
\endpmatrix.
$$

\vskip20pt

\centerline{Case $t=9$.}

\vskip10pt

The space $RJ_9$ has the bases
$$
\xi_{0,9}^{(1)}=-{\phi_{0,1}}({\phi_{0,4}})^2+
6{\phi_{0,2}}{\phi_{0,3}}{\phi_{0,4}}-5({\phi_{0,3}})^3=
$$
$$
\split
=&(3r[-1] - 2 + 3r^{-1}[-1]) +(-4r^6 + 6r^5 - 12r^4 +
22r^3 - 30r^2 + 36r -36 +\\
&36r^{-1}-30r^{-2}+ 22r^{-3} - 12r^{-4} + 6r^{-5} - 4r^{-6})q +\\
&(-r^9[-9] - 6r^8 + 15r^7 - 36r^6 + 72r^5 - 120r^4 + 171r^3 - 216r^2 +
255r -\\
&268 +255r^{-1} - 216r^{-2} + 171r^{-3} - 120r^{-4} +
72r^{-5} - 36r^{-6} + 15r^{-7} -\\
&6r^{-8} -r^{-9}[-9])q^2 + O(q^3);
\endsplit
$$
$$
\xi_{0,9}^{(2)}={\phi_{0,1}}({\phi_{0,4}})^2-
5{\phi_{0,2}}{\phi_{0,3}}{\phi_{0,4}}+4({\phi_{0,3}})^3=
$$
$$
\split
=&(r^2[-4]-r[-1]+4-r[-1]+r^{-2}[-4])+
(3r^6-8r^5+9r^4-24r^3+\\
&31r^2-32r+42-32r^{-1}+31r^{-2}-24r^{-3}+9r^{-4}-8r^{-5}+3r^{-6})q+\\
&(r^9[-9]+7r^8-15r^7+33r^6-80r^5+110r^4-177r^3+219r^2-241r+286-\\
&241r^{-1}+219r^{-2}-177r^{-3}+110r^{-4}-80r^{-5}+33r^{-6}-15r^{-7}+7r^{-8}+\\
&r^{-9}[-9])q^2+O(q^3);
\endsplit
$$
$$
\split 
\xi_{0,9}^{(3)}=&{E_{4,2}}{E_{4,3}}\Big({E_{4,1}}{\phi_{0,3}}-
{E_4}{\phi_{0,4}}\Big)/{\Delta_{12}}-
3{\phi_{0,1}}{\phi_{0,2}}({\phi_{0,3}})^2+\\
&2({\phi_{0,1}})^2{\phi_{0,3}}{\phi_{0,4}}-
30{\phi_{0,1}}({\phi_{0,4}})^2+
27{\phi_{0,2}}{\phi_{0,3}}{\phi_{0,4}}+9({\phi_{0,3}})^3=
\endsplit 
$$
$$
\split
=&q^{-1}[-36] + 24 +
(33r^6 + 486r^5 + 3159r^4 + 10758r^3 + 24057r^2 + 37908r +\\
&44082 +37908r^{-1} + 24057r^{-2} + 10758r^{-3} + 3159r^{-4} + 486r^{-5} +
33r^{-6})q +\\
&(2r^9[-9] + 243r^8 + 5346r^7 + 44055r^6 + 204120r^5 + 642978r^4 +
1483416r^3 +\\
&2632905r^2 + 3679020r + 4109590 +
3679020r^{-1} + 2632905r^{-2} +\\
&1483416r^{-3} + 642978r^{-4} +204120r^{-5} + 44055r^{-6} + 5346r^{-7} +\\ 
&243r^{-8} + 2r^{-9}[-9])q^2 + O(q^3).
\endsplit
$$
We have
$\overline{R}=\overline{P(\M_0)}$ where
$$
P(\M_0)=\pmatrix
{1}&{2}&{0}\cr
{0}&{-1}&{0}\cr
{-1}&{0}&{1}\cr
{2}&{9}&{1}\cr
\endpmatrix\
\equiv
\left[\matrix
4 & \overline{2}\\
1&   \overline{1}\\
36& \overline{0}\\
9&\overline{9}
\endmatrix\right];
\hskip10pt
Mul(P(\M_0),\,\xi)=
\pmatrix
0&1&0\\
3&0&0\\
0&0&1\\
-1&1&3
\endpmatrix;
$$
$$
G(P(\M_0))=
\pmatrix
{4}&{-2}&{-18}&{0}\cr
{-2}&{1}&{0}&{-9}\cr
{-18}&{0}&{36}&{-18}\cr
{0}&{-9}&{-18}&{9}\cr
\endpmatrix.
$$

\vskip20pt

\centerline{The case $t=12$.}

\vskip10pt

The space $RJ_{12}$ has the bases
$$
\xi_{0,12}^{(1)}=\bigl(\vartheta(\tau,z)/\eta(\tau)\bigr)^{12}=
$$
$$
\split
= &(r^8[-16] - 8r^7[-1] + 24r^6 - 24r^5 - 36r^4 +
120r^3 -\\
&88r^2 - 88r + 198 - 88r^{-1} - 88r^{-2} +\\
&120r^{-3} - 36r^{-4} - 24r^{-5} + 24r^{-6} - 8r^{-7}[-1] + r^{-8}[-16])q+\\
&(-4r^{10}[-4] + 24r^9 - 32r^8 - 104r^7 + 396r^6 - 352r^5 -512r^4 +
1440r^3 -\\
&904r^2 - 1008r + 2112 - 1008r^{-1} - 904r^{-2} + 1440r^{-3} -
512r^{-4} -\\
&352r^{-5} + 396r^{-6} - 104r^{-7} - 32r^{-8} + 24r^{-9} -
4r^{-10}[-4])q^2 + O(q^3)
\endsplit;
$$
$$
\xi_{0,12}^{(2)}=3{\phi_{0,2}}({\phi_{0,3}})^2{\phi_{0,4}}-
({\phi_{0,2}})^2({\phi_{0,4}})^2-2({\phi_{0,3}})^4-({\phi_{0,4}})^3=
$$
$$
\split
=&(r[-1] - 1 + r^{-1}[-1]) +
(-r^7[-1] + r^6 - r^5 + r^4 - r^2 +2r - 2 +\\
&2r^{-1} -r^{-2} + r^{-4} - r^{-5} + r^{-6} - r^{-7}[-1])q+
(-r^{10}[-4] + r^8 -
2r^7 + 3r^6 -\\
&3r^5 +2r^4 - 2r^2 +5r -6 + 5r^{-1} - 2r^{-2} + 2r^{-4} - 3r^{-5} +
3r^{-6} -\\
&2r^{-7} + r^{-8} - r^{-10}[-4])q + O(q^3)
\endsplit;
$$

$$
\xi_{0,12}^{(3)}=2({\phi_{0,2}})^2({\phi_{0,4}})^2-
5{\phi_{0,2}}({\phi_{0,3}})^2{\phi_{0,4}}+3({\phi_{0,3}})^4+({\phi_{0,4}})^3=
$$
$$
\split
= &(r^2[-4] - r[-1] + 3 - r^{-1}[-1] + r^{-2}[-4]) +
(r^7[-1] - 3r^6 + r^5 - 3r^4 +\\
&3r^3 -2r +6 -2r^{-1}+ 3r^{-2} - 3r^{-4} + r^{-5} - 3r^{-6} +
r^{-7}[-1])q +\\
&(2r^{10}[-4] - 3r^8 + 2r^7 - 9r^6 + 3r^5 - 6r^4 + 7r^2 - 5r +
18 -5r^{-1} +\\
&7r^{-2} - 6r^{-4} + 3r^{-5} - 9r^{-6} + 2r^{-7} - 3r^{-8} +
2r^{-10}[-4])q^2+O(q^3)
\endsplit;
$$

$$
\split 
\xi_{0,12}^{(4)}=&
{E_{4,3}}\Big({E_{4,1}}{E_{4,2}}({\phi_{0,3}})^2-
2{E_4}{E_{4,2}}{\phi_{0,3}}{\phi_{0,4}}+
{E_4}{E_{4,1}}({\phi_{0,4}})^2\Big)/{\Delta_{12}}-\\
&2({\phi_{0,1}})^2{\phi_{0,2}}({\phi_{0,4}})^2+
5({\phi_{0,1}})^2({\phi_{0,3}})^2{\phi_{0,4}}-
3{\phi_{0,1}}{\phi_{0,2}}({\phi_{0,3}})^3-\\
&36{\phi_{0,1}}{\phi_{0,3}}({\phi_{0,4}})^2+
24{\phi_{0,2}}({\phi_{0,3}})^2{\phi_{0,4}}+
9({\phi_{0,3}})^4+16({\phi_{0,4}})^3=
\endsplit 
$$
$$
\split
=& q^{-1}[-48] + 24 + (24r^7[-1] + 264r^6 +
1608r^5 + 5610r^4 +13464r^3 +\\
&24312r^2 + 34056r + 38208 + 34056r^{-1} + 24312r^{-2} + 13464r^{-3} +\\
&5610r^{-4} + 1608r^{-5} + 264r^{-6} + 24r^{-7}[-1])q +\\
&(12r^{10}[-4] + 440r^9 + 5544r^8 + 34104r^7 + 135388r^6 +
395808r^5 +\\
&902352r^4 + 1667360r^3 + 2550552r^2 + 3276240r + 3558160 +3276240r^{-1} +\\
&2550552r^{-2} + 1667360r^{-3} + 902352r^{-4} + 395808r^{-5} +
135388r^{-6} +\\
&34104r^{-7} + 5544r^{-8} + 440r^{-9} + 12r^{-10}[-4])q^2+O(q^3)
\endsplit.
$$
We have
$$
P(\M_0)=\pmatrix
{1}&{2}&{0}\cr
{0}&{-1}&{0}\cr
{-1}&{0}&{1}\cr
{1}&{8}&{1}\cr
\endpmatrix
\equiv
\left[\matrix
4&\overline{2}\\
1&\overline{1}\\
48&\overline{0}\\
16&\overline{8}
\endmatrix\right];\  
G(P(\M_0))=\pmatrix
{4}&{-2}&{-24}&{-8}\cr
{-2}&{1}&{0}&{-8}\cr
{-24}&{0}&{48}&{0}\cr
{-8}&{-8}&{0}&{16}\cr
\endpmatrix.
$$
$$
\overline{R}=
\left[\matrix
4,\overline{2}\\
1,\overline{1}\\
48,\overline{0}\\
16,\overline{8}\\
1,\overline{7}\\
4,\overline{10}
\endmatrix\right];
\ \ \
Mul(\overline{R},\,\xi)=
\pmatrix
0  & 0 & 1 & 0\\
0  & 1 & 0 & 0\\
0  & 0 & 0 & 1\\
1  & 0 & 0 & 0\\
-4 &-1 & 2 &12\\
-12&-2 & 3 &36
\endpmatrix;
$$

\vskip20pt

\centerline{Case $t=16$.}

\vskip10pt

The space $RJ_{16}$ has the bases
$$
\xi_{0,16}^{(1)}=\phi_{0,4}(\tau,2z)=
$$
$$
\split
=&(r^2[-4] +1 + r^{-2}[-4]) +
(-r^8 - r^6 + r^2 + 2 + r^{-2} - r^{-6} -\\
&r^{-8})q +(-r^{10} - 2r^8 - 2r^6 + 3r^2 + 4 + 3r^{-2} - 2r^{-6} - 2r^{-8} -
r^{-10})q^2 +\\
&(r^{14}[-4] -2r^{10} - 4r^8 - 4r^6 + 5r^2 + 8 +
5r^{-2} - 4r^{-6} - 4r^{-8} -\\
&2r^{-10} + r^{-14}[-4])q^3 + O(q^4)
\endsplit;
$$
$$
\split 
\xi_{0,16}^{(2)}=
&{E_{4,3}}\Big({E_4}{E_{4,1}}({\phi_{0,3}})^4-
({E_4})^2({\phi_{0,3}})^3{\phi_{0,4}}-\\
&2{E_{4,1}}{E_{4,2}}({\phi_{0,3}})^2{\phi_{0,4}}+
{E_4}{E_{4,2}}{\phi_{0,3}}({\phi_{0,4}})^2-
{E_{4,1}}{E_{4,2}}({\phi_{0,3}})^2{\phi_{0,4}} +\\
&2{E_4}{E_{4,2}}{\phi_{0,3}}({\phi_{0,4}})^2-
{E_4}{E_{4,1}}({\phi_{0,4}})^3\Big)/{\Delta_{12}}+\\
&2({\phi_{0,1}})^3({\phi_{0,3}})^3{\phi_{0,4}}-
3({\phi_{0,1}})^2{\phi_{0,2}}({\phi_{0,3}})^4-
7({\phi_{0,1}})^2({\phi_{0,3}})^2({\phi_{0,4}})^2-\\
&31{\phi_{0,1}}{\phi_{0,2}}({\phi_{0,3}})^3{\phi_{0,4}}+
46{\phi_{0,1}}({\phi_{0,3}})^5+
72{\phi_{0,1}}{\phi_{0,3}}({\phi_{0,4}})^3+\\
&7({\phi_{0,2}})^3({\phi_{0,3}})^2{\phi_{0,4}}-
72{\phi_{0,2}}({\phi_{0,3}})^2({\phi_{0,4}})^2-
197({\phi_{0,3}})^4{\phi_{0,4}}+\\
&2({\phi_{0,1}})^2{\phi_{0,2}}({\phi_{0,4}})^3+
21({\phi_{0,3}})^4{\phi_{0,4}}-26({\phi_{0,4}})^4+\\
&2{\phi_{0,1}}({\phi_{0,2}})^2{\phi_{0,3}}({\phi_{0,4}})^2-
({\phi_{0,2}})^2({\phi_{0,3}})^4-4({\phi_{0,2}})^2({\phi_{0,4}})^3-
2({\phi_{0,2}})^4({\phi_{0,4}})^2=
\endsplit 
$$
$$
\split
=&q^{-1}[-64]+((8r[-1] +14 + 8r^{-1}[-1]) +
(21r^8 + 200r^7 + 1036r^6 + 3360r^5 +\\
&8100r^4 + 15240r^3 + 23604r^2 + 30352r + 33058 +
30352r^{-1} + 23604r^{-2} +\\
&15240r^{-3} + 8100r^{-4} + 3360r^{-5} + 1036r^{-6} + 200r^{-7} +
21r^{-8})q +\\
&(56r^{11} + 1008r^{10} + 7336r^9 + 32932r^8 + 108800r^7 +
283504r^6 +610344r^5 +\\
&1112832r^4 + 1750728r^3 + 2401952r^2 + 2896688r^1 +
3081400 + \cdots)q^2+\\
&(4r^{14}[-4] +560r^{13} + 8092r^{12} + 58328r^{11} +
283784r^{10} +1042328r^9 +\\
&3082176r^8 + 7616904r^7 + 16136000r^6 + 29802144r^5 +
48582612r^4 +\\
&70497736r^3 +91619124r^2 + 107054192r +
112732002 + \cdots)q^3+O(q^4)
\endsplit.
$$
We have $\overline{R}=\overline{P(\M_0)}$ where
$$
P(\M_0)=
\pmatrix
{1}&{2}&{0}\cr
{0}&{-1}&{0}\cr
{-1}&{0}&{1}\cr
{5}&{32}&{3}\cr
{3}&{14}&{1}\cr
\endpmatrix
\equiv
\left[\matrix
4&\overline{2}\\
1&\overline{1}\\
64&\overline{0}\\
64&\overline{0}\\
4&\overline{14}\\
\endmatrix\right];
\hskip10pt
Mul(P(\M_0),\,\xi)=
\pmatrix
1 & 0\\
1&8\\
0&1\\
0&1\\
1&4
\endpmatrix;
$$
$$
G(P(\M_0))=
\pmatrix
{4}&{-2}&{-32}&{-32}&{-4}\cr
{-2}&{1}&{0}&{-32}&{-14}\cr
{-32}&{0}&{64}&{-64}&{-64}\cr
{-32}&{-32}&{-64}&{64}&{0}\cr
{-4}&{-14}&{-64}&{0}&{4}\cr
\endpmatrix.
$$

\vskip20pt

\centerline{Case $t=36$.}

\vskip10pt

The space $RJ_{36}$ has the bases
$$
\split
\xi_{0,36}^{(1)} =&
(-3r[-1] +5 - 3r^{-1}[-1]) +
(r^{12} + 3r^{11}+\cdots )q+\\
&(r^{18}[-36] - 3r^{17}[-1] +9r^{16} +\cdots)q^2+
(6r^{20} - 3r^{19}+\cdots )q^3+\\
&(4r^{24} - 15r^{22} +\cdots )q^4+\\
&((3r^{27}[-9] - 9r^{26}+3r^{25}+\cdots)q^5+
(3r^{29} + 6r^{28} +\cdots )q^6+\\
&(3r^{32}[-16] -25r^{30} + 9r^{29}+\cdots)q^7+
(-3r^{33} + 33r^{32}+\cdots )q^8+O(q^9)
\endsplit;
$$
$$
\split
\xi_{0,36}^{(2)}=&((r^2[-4] - r[-1] + 1 - r[-1] + r^2[-4]) +
(-r^{12} + r^{11} - r^{10} + \cdots )q+\\
&(-r^{17}[-1] + r^{16} - r^{15}+\cdots )q^2+
(-r^{19} + 2r^{18} - 3r^{17} +\cdots )q^3+\\
&(-r^{21}  + 2r^{20} - 4r^{19}+\cdots )q^4+\\
&((r^{27}[-9] -r^{26} + r^{25}+\cdots )q^5+
(r^{29} - 2r^{28}  + 3r^{27} + \cdots )q^6+\\
&(r^{32}[-16] -r^{30}+ 3r^{29}+\cdots )q^7+\\
&(r^{34}[-4] - r^{33} + r^{32} - 3r^{30}+\cdots)q^8 + O(q^9)
\endsplit;
$$
$$
\split
\xi_{0,36}^{(3)}=&q^{-1}[-144] + 24 +
(24r^{12} + 72r^{11} +\cdots)q+\\
&4r^{18}[-36] + 144r^{16} + 672r^{15}+\cdots)q^2+
(144r^{20} + 1008r^{19}+\cdots)q^3+\\
&(24r^{24}+ 288r^{23}+\cdots)q^4+\\
&(8r^{27}[-9] +216r^{26} + 3096r^{25}+\cdots )q^5+\\
&(72r^{29} + 1584r^{28} + 15720r^{27}+\cdots)q^6+\\
&(9r^{32}[-16] + 288r^{31} + 5304r^{30}+\cdots)q^7+\\
&(672r^{33}+ 12096r^{32}+\cdots)q^8+O(q^9)
\endsplit.
$$
We have
$$
P(\M_0)=
\left(\smallmatrix
{1}&{2}&{0}\cr
{0}&{-1}&{0}\cr
{-1}&{0}&{1}\cr
{2}&{18}&{1}\cr
{5}&{27}&{1}\cr
{7}&{32}&{1}\cr
\endsmallmatrix\right);\ \ \
G(P(\M_0))=\left(\smallmatrix
{4}&{-2}&{-72}&{-36}&{-18}&{-8}\cr
{-2}&{1}&{0}&{-18}&{-27}&{-32}\cr
{-72}&{0}&{144}&{-72}&{-288}&{-432}\cr
{-36}&{-18}&{-72}&{36}&{-18}&{-72}\cr
{-18}&{-27}&{-288}&{-18}&{9}&{0}\cr
{-8}&{-32}&{-432}&{-72}&{0}&{16}\cr
\endsmallmatrix\right);
$$
$$
\overline{R}=
\left[\matrix
1,\ \overline{1}\\
1,\ \overline{17}\\
4,\ \overline{2}\\
4,\ \overline{34}\\
9,\ \overline{27}\\
16, \overline{32}\\
36, \overline{18}\\
144, \overline{0}
\endmatrix\right];
\hskip10pt
Mul(\overline{R},\,\xi)=
\pmatrix
-3  &0  & 0\\
-3  & 0 & 0\\
 0  & 1 & 0\\
 0  & 1 & 0\\
 4  & 1 & 12\\
 3  & 1 & 9\\
 1  & 0 & 4\\
 0  & 0 & 1
\endpmatrix .
$$

\head
6. The list of all Lorentzian Kac--Moody algebras with the
denominator identity function which is automorphic with respect to
the extended paramodular group $\widehat{O}(L_t)$.
\endhead

In Table 2 below we give the list of Lorentzian Kac--Moody algebras from
Theorem 1.1. The product part of their denominator identities is
defined by the infinite products of
Theorem 2.1 for some Jacobi forms $\xi$ from the spaces
$RJ_t$ which were described in Table 1 by their basis.
These products are characterised by the property that
multiplicities of divisors of the infinite product
$B_\xi$ are equal to 0 or 1
for any rational quadratic divisor which is orthogonal to
a root of $L_t$. Since $B_\xi$ is reflective, it is the whole
divisor of $B_\xi$. We denote the corresponding Lorentzian
Kac--Moody algebra as $\geg(\xi)$ since it is defined by
the Jacobi form $\xi$.

We describe the fundamental polyhedron $\M$ and
the set $P(\M)$ of orthogonal roots to $\M$ defining the Weyl
group and the set of simple real roots of the Algebra
$\geg(\xi)$. We also give the
subset $P(\M)_{\1o}\subset P(\M)$ of super roots. 
It means that the corresponding 
generators $e_\alpha$, $f_\alpha$, $\alpha \in P(\M)_{\1o}$, 
should be super.
If we don't mention the set $P(\M)_{\1o}$, it is empty.
We also give the generalised Cartan matrix
$$
A=\left({2(\alpha_i,\alpha_j)\over \alpha_i^2}\right),\ \ 
\alpha_i,\,\alpha_j \in P(\M), 
$$
which is the main invariant of the algebra.
We also give the Weyl vector $\rho$.

All these polyhedra $\M$ are composed from the fundamental polyhedron
$\M_0$ for $W(S_t)$ using some group of symmetries of the polyhedron $\M$.
We use these symmetries to describe the sets $P(\M)$ and $P(\M)_{\1o}$
using the set $P(\M_0)$. We numerate as $\alpha_1,\dots ,\alpha_k$
the elements of $P(\M)$ as they are given in Table 1.
We denote by $s_\alpha$ the reflection in the root $\alpha$. It is
given by the formula
$$
s_\alpha\,:\,x\to x-{2(x,\,\alpha )\over \alpha^2}\alpha,\ \  
x\in S_t^\ast.
$$
We denote by $[g_1,\dots ,g_k]$ the group generated by $g_1,\dots ,g_k$

\vskip40pt

\centerline{\bf Table 2. The list of all 29 Lorentzian Kac--Moody algebras}
\centerline{\bf of the rank three with the root lattice $S_t^\ast$ and}
\centerline{\bf the symmetry group $\widehat{O}^+(L_t)$ (from Theorem 1.1).} 

\vskip20pt

\centerline{Case $t=1$}.

\vskip10pt

The Algebra $\geg(\xi_{0,1}^{(1)})$.
The fundamental chamber
$\M=[s_{\alpha_1},\,s_{\alpha_3}](\M_0)$
is the right triangle with zero angles. We have
$$
P(\M)=[s_{\alpha_1},\,s_{\alpha_3}](\alpha_2)
$$
with the group of symmetries $[s_{\alpha_1},s_{\alpha_3}]$ which is
$D_3$. The generalised Cartan matrix is
$$
A_{1,II}=
\pmatrix
\hphantom{-}{2}&{-2}&{-2}\cr
{-2}&\hphantom{-}{2}&{-2}\cr
{-2}&{-2}&\hphantom{-}{2}\cr
\endpmatrix.
$$
The Weyl vector
$\rho=({1\over 2},{1\over 2},{1\over 2})$.
The corresponding automorphic form $B_{\xi_{0,1}^{(1)}}$
is $\Delta_5$ of the weight 5 which is product of 
ten even theta-constants of genus 2. See \cite{GN1}, \cite{GN2}, \cite{GN5}.

The Algebra $\geg(\xi_{0,1}^{(2)})$.
The chamber $\M=\M_0$ is a triangle with
angles $\pi/3$, $0$, $\pi/2$. The set
$P(\M)=P(\M_0)$, $P(\M)_{\1o}=\{\alpha_2\}$. The generalised
Cartan matrix is
$$
A_{1,I,\0o}=
\pmatrix
2 & -1 & -1\\
-4&  2 &  0\\
-1&  0 &  2
\endpmatrix .
$$
The Weyl vector $\rho=({5\over 2},{1\over 2},{3\over 2})$.
The corresponding automorphic form is  
Igusa's \cite{Ig}  
modular form $\Delta_{30}=\Delta_{35}/\Delta_5$ of the weight 30. 
See \cite{GN4}, \cite{GN5}.

The Algebra $\geg(\xi_{0,1}^{(1)}+\xi_{0,1}^{(2)})$.
The chamber $\M=\M_0$ and
$$
P(\M)=\{\alpha_1,\,2\alpha_2,\,\alpha_3\}.
$$
The generalised Cartan matrix is
$$
A_{1,0}=
\pmatrix
2&-2&-1\\
-2&2&0\\
-1&0&2
\endpmatrix.
$$
The Weyl vector $\rho=(3,\,1,\,2)$. The automorphic form
is Igusa's \cite{Ig} modular form $\Delta_{35}$ of the weight 35. 
See \cite{GN4} (for a new simple construction of this form) 
and \cite{GN5}.

\vskip20pt

\centerline{Case $t=2$.}

\vskip10pt

The Algebra $\geg(\xi_{0,2}^{(1)})$. The chamber
$\M=[s_{\alpha_1},\,s_{\alpha_3}](\M_0)$ is
the right quadrangle with zero angles;
$$
P(\M)=[s_{\alpha_1},\,s_{\alpha_3}](\alpha_2)
$$
with the group of symmetries $[s_{\alpha_1},\,s_{\alpha_3}]$
which is $D_4$. The generalised Cartan matrix is
$$
A_{2,II}=
\pmatrix
\hphantom{-}{2}&{-2}&{-6}&{-2}\cr
{-2}&\hphantom{-}{2}&{-2}&{-6}\cr
{-6}&{-2}&\hphantom{-}{2}&{-2}\cr
{-2}&{-6}&{-2}&\hphantom{-}{2}\cr
\endpmatrix.
$$
The Weyl vector $\rho=({1\over 4},\,{1\over 2},\,{1\over 4})$.
The corresponding automorphic form is $\Delta_2$ of the weight 2. 
See \cite{GN1} and \cite{GN5}.

The Algebra $\geg(\xi_{0,2}^{(2)})$.
The chamber $\M=[s_{\alpha_3}](\M_0)$ is a triangle with
angles $0,\,0,\,\pi/2$. The sets
$$
P(\M)=\{\alpha_1,\,\alpha_2,\,s_{\alpha_3}(\alpha_1)\},\ \ \
P(\M)_{\1o}=\{\alpha_2\}
$$
with the group of symmetries $[s_{\alpha_3}]$ which is $D_1$.
The generalised Cartan matrix is
$$
A_{2,I,\0o}=
\pmatrix
2 & -1 &  0\\
-4&  2 & -4\\
0 & -1 &  2
\endpmatrix. 
$$
The Weyl vector $\rho=({3\over 4},{1\over 2},{3\over 4})$.
The corresponding automorphic form is
$\Delta_{9}=\Delta_{11}/\Delta_2$ of the weight 9.
See \cite{GN4}, \cite{GN5}.

The Algebra $\geg(\xi_{0,2}^{(1)}+\xi_{0,2}^{(2)})$. The
polygon $\M=[s_{\alpha_3}](\M_0)$
is a triangle with angles $0,\,0,\,\pi/2$
(the same as for $\xi_{0,2}^{(2)}$); the set
$$
P(\M)=\{\alpha_1,\,2\alpha_2,\,s_{\alpha_3}(\alpha_1)\}
$$
with the group of symmetries $[s_{\alpha_3}]$ which is $D_1$.
The generalised Cartan matrix is
$$
A_{2,0}=
\pmatrix
2 & -2& 0\\
-2&  2&-2\\
 0&	-2&	2
\endpmatrix.
$$
The Weyl vector $\rho=(1,1,1)$. The corresponding automorphic form
is $\Delta_{11}$ of the weight 11. See \cite{GN4}, \cite{GN5}.

The Algebra $\geg(\xi_{0,2}^{(3)})$.
The chamber $\M=[s_{\alpha_1},\,s_{\alpha_2}](\M_0)$ is an
infinite polygon with angles $\pi/2$ and
which is touching a horosphere with the
centre at the Weyl vector $\rho=(1,0,0)$. The set
$$
P(\M)=[s_{\alpha_1},\,s_{\alpha_2}](\alpha_3)
$$
with the group of symmetries $[s_{\alpha_1},\,s_{\alpha_2}]$
which is $D_\infty$. The generalised
Cartan matrix is a symmetric matrix
$$
A_{2,\1o}=
\left({(\alpha,\,\alpha^\prime)\over 4}\right),\ \
\alpha,\alpha^\prime \in P(\M).
$$
The corresponding automorphic form is $\Psi_{12}^{(2)}$ of the 
weight 12.
See \cite{GN5}.

The Algebra $\geg(\xi_{0,2}^{(1)}+\xi_{0,2}^{(3)})$.
The polygon $\M=[s_{\alpha_1}](\M_0)$ is a quadrangle with angles
$0,\pi/2,\pi/2,\,\pi/2$; the set
$$
P(\M)=[s_{\alpha_1}]\{\alpha_2,\,\alpha_3\}=
\{\alpha_2,\alpha_3,(1,4,1),(1,1,0)\}
$$
with the group of symmetries $[s_{\alpha_1}]$ which is $D_1$.
The generalised Cartan matrix is
$$
A_{2,II,\1o}=
\pmatrix
{2}&{0}&{-8}&{-2}\cr
{0}&{2}&{0}&{-1}\cr
{-1}&{0}&{2}&{0}\cr
{-2}&{-8}&{0}&{2}\cr
\endpmatrix.
$$
The Weyl vector
$\rho=({{5}\over{4}},{{1}\over{2}},{{1}\over{4}})$.
The automorphic form is
$\Delta_{14}=\Delta_2\cdot \Psi_{12}^{(2)}$ of the weight 14.
(We must correct the case $(2,II,\1o)$ in \cite{GN5, page 264}
in this way.)

The Algebra $\geg(\xi_{0,2}^{(2)}+\xi_{0,2}^{(3)})$.
The polygon $\M=\M_0$ is the triangle with angles
$0,\,\pi/2,\pi/4$; the set
$$
P(\M)=P(\M_0),\ \ \  P(\M)_{\1o}=\{\alpha_2\}
$$
with the trivial group of symmetries and
with the generalised Cartan matrix
$$
A_{2,I,\1o}=
\pmatrix
2 & -1 & -2\\
-4&  2 &  0\\
-1&  0 &  2
\endpmatrix.
$$
The Weyl vector $\rho=({{7}\over{4}},{{1}\over{2}},{{3}\over{4}})$.
The automorphic form is $\Delta_9\cdot \Psi_{12}^{(2)}$ of the weight $21$.
See \cite{GN5}.

The Algebra $\geg(\xi_{0,2}^{(1)}+\xi_{0,2}^{(2)}+\xi_{0,2}^{(3)})$.
The polygon $\M=\M_0$ is the triangle with angles
$0,\pi/2,\,\pi/4$ (the same as for the
$\geg(\xi_{0,2}^{(2)}+\xi_{0,2}^{(3)})$); the set
$$
P(\M)=\{\alpha_1\,2\alpha_2,\alpha_3\}.
$$
with the trivial group of symmetries and
with the generalised Cartan matrix
$$
A_{2,0,\1o}=
\pmatrix
2 & -2 & -2\\
-2&  2 &  0\\
-1&  0 &  2
\endpmatrix.
$$
The Weyl vector $\rho=(2,\,1,\,1)$.
The automorphic form is
$\Delta_2\cdot \Delta_9\cdot \Psi_{12}^{(2)}$ of the weight $23$.
See \cite{GN5}.

\centerline{Case $t=3$.}

\vskip10pt

The Algebra $\geg(\xi_{0,3}^{(1)})$. The chamber
$\M=[s_{\alpha_1},\,s_{\alpha_3}](\M_0)$ is
the right hexagon with zero angles, 
$$
P(\M)=[s_{\alpha_1},\,s_{\alpha_3}](\alpha_2)
$$
with the group of symmetries $[s_{\alpha_1},\,s_{\alpha_3}]$
which is $D_6$. The generalised Cartan matrix is
$$
A_{3,II}=
\pmatrix
\hphantom{-}{2}&{-2}&{-10}&{-14}&{-10}&{-2}\cr
{-2}&\hphantom{-}{2}&{-2}&{-10}&{-14}&{-10}\cr
{-10}&{-2}&\hphantom{-}{2}&{-2}&{-10}&{-14}\cr
{-14}&{-10}&{-2}&\hphantom{-}{2}&{-2}&{-10}\cr
{-10}&{-14}&{-10}&{-2}&\hphantom{-}{2}&{-2}\cr
{-2}&{-10}&{-14}&{-10}&{-2}&\hphantom{-}{2}\cr
\endpmatrix . 
$$
The Weyl vector $\rho=({1\over 6},\,{1\over 2},\,{1\over 6})$.
The corresponding automorphic form is $\Delta_1$ of the weight 1. 
See \cite{GN4}, \cite{GN5}.

The Algebra $\geg(\xi_{0,3}^{(2)})$.
The chamber $\M=[s_{\alpha_3}](\M_0)$ is a triangle with
angles $0,\,0,\,\pi/3$. The sets
$$
P(\M)=\{\alpha_1,\,\alpha_2,\,s_{\alpha_3}(\alpha_1)\},\ \ \
P(\M)_{\1o}=\{\alpha_2\}
$$
with the group of symmetries $[s_{\alpha_3}]$ which is $D_1$.
The generalised Cartan matrix is
$$
A_{3,I,\0o}=
\pmatrix
2 & -1 & -1\\
-4&  2 & -4\\
-1& -1 &  2
\endpmatrix . 
$$
The Weyl vector $\rho=({1\over 2},{1\over 2},{1\over 2})$.
The corresponding automorphic form is
$D_6$ of the weight 6.
See \cite{GN4}, \cite{GN5}.

The Algebra $\geg(\xi_{0,3}^{(1)}+\xi_{0,3}^{(2)})$. The
polygon $\M=[s_{\alpha_3}](\M_0)$
is a triangle with angles $0,\,0,\,\pi/3$
(the same as for $\xi_{0,3}^{(2)}$); the set
$$
P(\M)=\{\alpha_1,\,2\alpha_2,\,s_{\alpha_3}(\alpha_1)\}
$$
with the group of symmetries $[s_{\alpha_3}]$ which is $D_2$.
The generalised Cartan matrix is
$$
A_{3,0}=
\pmatrix
2 & -2& -1\\
-2&  2&-2\\
-1&	-2&	2
\endpmatrix.
$$
The Weyl vector $\rho=({2\over 3},1,{2\over 3})$.
The corresponding automorphic form is $\Delta_1\cdot D_6$ of 
the weight 7. See \cite{GN4}, \cite{GN5}.

The Algebra $\geg(\xi_{0,3}^{(3)})$.
The chamber $\M=[s_{\alpha_1},\,s_{\alpha_2}](\M_0)$ is an
infinite polygon with angles $\pi/3$ and
which is touching a horosphere with the
centre at the Weyl vector $\rho=(1,0,0)$. The set
$$
P(\M)=[s_{\alpha_1},\,s_{\alpha_2}](\alpha_3)
$$
with the group of symmetries $[s_{\alpha_1},\,s_{\alpha_2}]$
which is $D_\infty$. The generalised
Cartan matrix is a symmetric matrix
$$
A_{3,\1o}=
\left({(\alpha,\,\alpha^\prime)\over 6}\right),\ \
\alpha,\alpha^\prime \in P(\M).
$$
The corresponding automorphic form is $\Psi_{12}^{(3)}$ of the 
weight 12. See \cite{GN5}.

The Algebra $\geg(\xi_{0,3}^{(1)}+\xi_{0,3}^{(3)})$.
The polygon $\M=[s_{\alpha_1}](\M_0)$ is a quadrangle with angles
$0,\pi/2,\pi/3,\,\pi/2$; the set
$$
P(\M)=[s_{\alpha_1}]\{\alpha_2,\,\alpha_3\}=
\{\alpha_2,\alpha_3,(2,6,1),(1,1,0)\}
$$
with the group of symmetries $[s_{\alpha_1}]$ which is $D_1$.
The generalised Cartan matrix is
$$
A_{3,II,\1o}=
\pmatrix
{2}&{0}&{-12}&{-2}\cr
{0}&{2}&{-1}&{-1}\cr
{-1}&{-1}&{2}&{0}\cr
{-2}&{-12}&{0}&{2}\cr
\endpmatrix . 
$$
The Weyl vector
$\rho=({{7}\over{6}},{{1}\over{2}},{{1}\over{6}})$.
The automorphic form is
$\Delta_1\cdot \Psi_{12}^{(3)}$ of the weight 13.
(We must correct the case $(3,II,\1o)$ in \cite{GN5, page 264}
in this way.)

The Algebra $\geg(\xi_{0,3}^{(2)}+\xi_{0,3}^{(3)})$.
The polygon $\M=\M_0$ is the triangle with angles
$0,\,\pi/2,\pi/6$; the set
$$
P(\M)=P(\M_0),\ \ \  P(\M)_{\1o}=\{\alpha_2\}
$$
with the trivial group of symmetries and
with the generalised Cartan matrix
$$
A_{3,I,\1o}=
\pmatrix
2 & -1 & -3\\
-4&  2 &  0\\
-1&  0 &  2
\endpmatrix.
$$
The Weyl vector $\rho=({{3}\over{2}},{{1}\over{2}},{{1}\over{2}})$.
The automorphic form is $D_6\cdot \Psi_{12}^{(3)}$ of the weight $18$.
See \cite{GN5}.

The Algebra $\geg(\xi_{0,3}^{(1)}+\xi_{0,3}^{(2)}+\xi_{0,3}^{(3)})$.
The polygon $\M=\M_0$ is the triangle with angles
$0,\pi/2,\,\pi/6$ (the same as for the
$\geg(\xi_{0,3}^{(2)}+\xi_{0,3}^{(3)})$); the set
$$
P(\M)=\{\alpha_1\,2\alpha_2,\alpha_3\}.
$$
with the trivial group of symmetries and
with the generalised Cartan matrix
$$
A_{3,0,\1o}=
\pmatrix
2 & -2 & -3\\
-2&  2 &  0\\
-1&  0 &  2
\endpmatrix.
$$
The Weyl vector $\rho=({5\over 3},\,1,\,{2\over 3})$.
The automorphic form is
$\Delta_1\cdot D_6\cdot \Psi_{12}^{(3)}$ of the weight $19$.
See \cite{GN5}.

\vskip20pt

\centerline{Case $t=4$.}

\vskip10pt

The Algebra $\geg(\xi_{0,4}^{(1)})$. The chamber
$\M=[s_{\alpha_1},\,s_{\alpha_3}](\M_0)$ is
the infinite polygon with zero angles touching a horosphere with the
centre at the Weyl vector $\rho=({1\over 8},{1\over 2},{1\over 8})$;
$$
P(\M)=[s_{\alpha_1},\,s_{\alpha_3}](\alpha_2)]
$$
with the group of symmetries $[s_{\alpha_1},\,s_{\alpha_3}]$
which is $D_\infty$. The generalised Cartan matrix is
$$
A_{4,II,\0o}=\left(2(\alpha,\alpha^\prime)\right),\ \ \
\alpha,\,\alpha^\prime \in P(\M).
$$
The corresponding automorphic form is $\Delta_{1/2}$ of the 
weight 1/2 which is the theta-constant of the genus 2. See \cite{GN5}. 

The Algebra $\geg(\xi_{0,4}^{(2)})$.
The chamber $\M=[s_{\alpha_3}](\M_0)$ is a triangle with
angles $0,\,0,\,0$. The sets
$$
P(\M)=\{\alpha_1,\,\alpha_2,\,s_{\alpha_3}(\alpha_1)\};\ \ \
P(\M)_{\1o}=\{\alpha_2\}
$$
with the group of symmetries $[s_{\alpha_3}]$ which is $D_1$.
The generalised Cartan matrix is
$$
A_{4,I,\0o}=
\pmatrix
2 & -1 &  -2\\
-4&  2 & -4\\
-2 & -1 &  2
\endpmatrix
$$
The Weyl vector $\rho=({3\over 8},{1\over 2},{3\over 8})$.
The corresponding automorphic form is
$\Delta_5^{(4)}/\Delta_{1/2}$ of the weight ${9\over 2}$.
See \cite{GN5}.

The Algebra $\geg(\xi_{0,4}^{(1)}+\xi_{0,4}^{(2)})$. The
polygon $\M=[s_{\alpha_3}](\M_0)$
is a triangle with angles $0,\,0,\,0$
(the same as for $\xi_{0,1}^{(1)}$ and $\xi_{0,4}^{(2)}$); the set
$$
P(\M)=\{\alpha_1,\,2\alpha_2,\,s_{\alpha_3}(\alpha_1)\}
$$
with the group of symmetries $[s_{\alpha_3}]$ which is $D_1$.
The generalised Cartan matrix is
$$
A_{4,0,\0o}=A_{1,II}=
\pmatrix
2 & -2& -2\\
-2&  2&-2\\
-2& -2& 2
\endpmatrix
$$
(it is the same as for $\xi_{0,1}^{(1)}$). 
The Weyl vector $\rho=({1\over 2},1,{1\over 2})$.
The corresponding automorphic form
is 
$\Delta_5^{(4)}(z_1,z_2,z_3)=\Delta_5(z_1,2z_2,z_3)$ 
of the weight 5. This case is equivalent to 
the case $\geg(\xi_{0,1}^{(1)})$ above. See \cite{GN5}. 

The Algebra $\geg(\xi_{0,4}^{(3)})$.
The chamber $\M=[s_{\alpha_1},\,s_{\alpha_2}](\M_0)$ is an
infinite polygon with zero angles
which is touching a horosphere with the
centre at the Weyl vector $\rho=(1,0,0)$. The set
$$
P(\M)=[s_{\alpha_1},\,s_{\alpha_2}](\alpha_3)
$$
with the group of symmetries $[s_{\alpha_1},\,s_{\alpha_2}]$
which is $D_\infty$. The generalised
Cartan matrix is a symmetric matrix
$$
A_{4,\1o}=
\left({(\alpha,\,\alpha^\prime)\over 8}\right),\ \
\alpha,\alpha^\prime \in P(\M).
$$
The corresponding automorphic form is $\Psi_{12}^{(4)}$ of the 
weight 12. See \cite{GN5}.

The Algebra $\geg(\xi_{0,4}^{(1)}+\xi_{0,4}^{(3)})$.
The polygon $\M=[s_{\alpha_1}](\M_0)$ is a quadrangle with angles
$0,\pi/2,\,0,\,\pi/2$; the set
$$
P(\M)=[s_{\alpha_1}]\{\alpha_2,\,\alpha_3\}=
\{\alpha_2,\alpha_3,(3,8,1),(1,1,0)\}
$$
with the group of symmetries $[s_{\alpha_1}]$ which is $D_1$.
The generalised Cartan matrix is
$$
A_{4,II,\1o}=
\pmatrix
{2}&{0}&{-16}&{-2}\cr
{0}&{2}&{-2}&{-1}\cr
{-1}&{-2}&{2}&{0}\cr
{-2}&{-16}&{0}&{2}\cr
\endpmatrix .
$$
The Weyl vector is
$({9\over 8},{1\over 2},{1\over 8})$.
The automorphic form is
$\Delta_{1/2}\cdot \Psi_{12}^{(4)}$ of the weight ${25\over 2}$.
(We must correct the case $(2,II,\1o)$ in \cite{GN5, page 264}
in this way.)

The Algebra $\geg(\xi_{0,4}^{(2)}+\xi_{0,4}^{(3)})$.
The polygon $\M=\M_0$ is the triangle with angles
$0,\,\pi/2,\,0$; the set
$$
P(\M)=P(\M_0),\ \ \  P(\M)_{\1o}=\{\alpha_2\}
$$
with the trivial group of symmetries and
with the generalised Cartan matrix
$$
A_{4,I,\1o}=
\pmatrix
2 & -1 & -4\\
-4&  2 &  0\\
-1&  0 &  2
\endpmatrix.
$$
The Weyl vector $\rho=({{11}\over{8}},{{1}\over{2}},{{3}\over{8}})$.
The automorphic form is
$\Psi_{12}^{(4)}\cdot \Delta_5^{(4)}/\Delta_{1/2}$ of the weight
${33\over 2}$.
See \cite{GN5}.

The Algebra $\geg(\xi_{0,4}^{(1)}+\xi_{0,4}^{(2)}+\xi_{0,4}^{(3)})$.
The polygon $\M=\M_0$ is the triangle with angles
$0,\pi/2,\,\pi/4$ (the same as for the
$\geg(\xi_{0,4}^{(2)}+\xi_{0,4}^{(3)})$); the set
$$
P(\M)=\{\alpha_1\,2\alpha_2,\alpha_3\}.
$$
with the trivial group of symmetries and
with the generalised Cartan matrix
$$
A_{4,0,\1o}=
\pmatrix
2 & -2 & -4\\
-2&  2 &  0\\
-1&  0 &  2
\endpmatrix.
$$
The Weyl vector $\rho=({3\over 2},\,1,\,{1\over 2})$.
The automorphic form is
$\Delta_5^{(4)}\cdot \Psi_{12}^{(4)}$ of the weight $17$.
See \cite{GN5}.

\vskip20pt

\centerline{Case $t=8$.}

\vskip10pt

The Algebra $\geg(\xi_{0,8}^{(2)})$. The chamber
$\M=[s_{\alpha_3}](\M_0)$ is the right quadrangle with
zero angles; the set
$$
P(\M)=[s_{\alpha_3}](\alpha_1,\,2\alpha_2,\,\alpha_4)=
\{\alpha_1,\,2\alpha_2,\,s_{\alpha_3}(\alpha_1),\,\alpha_4 \}
$$
with the group of symmetries $[s_{\alpha_3}]$ which is $D_1$.
The generalised Cartan matrix is
$A_{2,II}$ (the same as for $\geg(\xi_{0,2}^{(1)})$ for $t=2$).
The Weyl vector $\rho=({1\over 4},\,1,\,{1\over 4})$.
The automorphic form is
$\Delta_2^{(8)}(z_1,z_2,z_3)=\Delta_2(z_1,2z_2,z_3)$ of the 
weight 2 where $\Delta_2$ corresponds to $\geg(\xi_{0,2}^{(1)})$. 
This case is equivalent to $\geg(\xi_{0,2}^{(1)})$. 

\vskip20pt

\centerline{Case $t=9$.}

\vskip10pt

The algebra $\geg(\xi_{0,9}^{(2)})$. The chamber
$\M=[s_{\alpha_3}](\M_0)$ is the pentagon with angles
$0$, $0$, $\pi/2$, $0$, $\pi/2$; the set
$$
P(\M)=[s_{\alpha_3}](\alpha_1,\alpha_2,\alpha_4)=
\{\alpha_1,\alpha_2,s_{\alpha_3}(\alpha_1),s_{\alpha_3}(\alpha_4),
\alpha_4\},\ \ P(\M)_{\1o}=\{\alpha_2\}, 
$$
with the group of symmetries
$[s_{\alpha_2}]$ which is $D_1$.
The generalised Cartan matrix is
$$
\pmatrix
{2}&{-1}&{-7}&{-9}&{0}\cr
{-4}&{2}&{-4}&{-18}&{-18}\cr
{-7}&{-1}&{2}&{0}&{-9}\cr
{-4}&{-2}&{0}&{2}&{-2}\cr
{0}&{-2}&{-4}&{-2}&{2}\cr
\endpmatrix .
$$
The Weyl vector is $({1\over 6},{1\over 2},{1\over 6})$.
The automorphic form is $D_2$ of the weight 2. See \cite{GN5}.

\vskip20pt

\centerline{Case $t=12$.}

\vskip10pt

The Algebra $\geg(\xi_{0,12}^{(2)}+\xi_{0,12}^{(3)})$.
The chamber $\M=[s_{\alpha_3},s_{\alpha_4}](\M_0)$ is the
right hexagon with zero angles; the set
$$
\split
P(\M)=&[s_{\alpha_3},s_{\alpha_4}](\alpha_1,\,2\alpha_2)=\\
&\{\alpha_1,2\alpha_2,s_{\alpha_3}(\alpha_1),
s_{\alpha_4}s_{\alpha_3}(\alpha_1),s_{\alpha_4}(2\alpha_2),
s_{\alpha_4}(\alpha_1)\}
\endsplit
$$
with the group of symmetries $[s_{\alpha_3},\,s_{\alpha_4}]$ which
is $D_4$. The generalised Cartan matrix is
$A_{3,II}$ (the same as for $\geg(\xi_{0,3}^{(1)})$).
The Weyl vector is $({1\over 6},1,{1\over 6})$.
The automorphic form is $\Delta_1^{(12)}(z_1,z_2,z_3)=
\Delta_1(z_1,2z_2,z_3)$ of the weight 1 where $\Delta_1$ 
corresponds to $\geg(\xi_{0,3}^{(1)})$.  
This case is equivalent to $\geg(\xi_{0,3}^{(1)})$  
above.

\vskip20pt

\centerline{Case $t=16$.}

\vskip10pt

The Algebra $\geg(\xi_{0,16}^{(1)})$.
The chamber is the infinite polygon
$\M=[s_{\alpha_3},s_{\alpha_4}](\M_0)$ with zero angles
touching a horosphere with the centre at the Weyl vector
$\rho=({1\over 8},1,{1\over 8})$. The set
$$
P(\M)=[s_{\alpha_3},s_{\alpha_4}](\alpha_1,2\alpha_2,\alpha_5)
$$
with the group of symmetry $[s_{\alpha_3},s_{\alpha_4}]$ which is
$D_\infty$. The generalised Cartan matrix is
$$
\left({(\alpha,\,\alpha^\prime)\over 2}\right),\ \ \
\alpha,\ \alpha^\prime \in P(\M), 
$$
which is the same as for $\geg(\xi_{0,4}^{(1)})$.
The automorphic form is $\Delta_{1/2}^{(16)}(z_1,z_2,z_3)=
\Delta_{1/2}(z_1,2z_2,z_3)$ of the weight 1/2 
where $\Delta_{1/2}$ corresponds to
$\geg(\xi_{0,4}^{(1)})$. This case is equivalent to $\geg(\xi_{0,4}^{(1)})$.

\vskip20pt

\centerline{Case $t=36$.}

\vskip10pt

The Algebra $\geg(\xi_{0,36}^{(2)})$. The chamber
$\M=[s_{\alpha_3},s_{\alpha_4}](\M_0)$  
is the infinite periodic polygon  
with angles 
$\dots$, $0$, $\pi/2$, $0$, $0$, $0$, $0$, $\pi/2$, $0$, $\dots$, 
with the centre at the Weyl vector 
$\rho=({1\over 24},{1\over 2},{1\over 24})$ at infinity. 
The set 
$$
P(\M)=[s_{\alpha_3},\,s_{\alpha_4}]
(\alpha_1,\,\alpha_2,\,\alpha_5,\,\alpha_6),\ \ \
P(\M)_{\1o}=[s_{\alpha_3},\,s_{\alpha_4}](\alpha_2)
$$
with the group of symmetry $[s_{\alpha_3},\,s_{\alpha_4}]$ which is 
$D_\infty$. The generalised Cartan matrix is
$$
\left({2(\alpha,\alpha^\prime)\over (\alpha,\,\alpha)}\right),
\ \ \ \alpha,\ \alpha^\prime \in P(\M).
$$
The automorphic form is $D_{1/2}$ of the weight 1/2. See \cite{GN5}.

\head
7. Possible Physical Applications
\endhead 

It would be interesting to find some Quantum Systems with  
symmetries corresponding to the Lorentzian Kac--Moody 
algebras of Theorem 1.1 and to possible algebras with 
denominator identities which are reflective automorphic forms 
of Main Theorem 2.2. See \cite{B3}, \cite{DMVV}, \cite{DVV}, 
\cite{G2}, \cite{G3}, \cite{HM}, \cite{GN3}, \cite{GN6}, 
\cite{Ka}, \cite{KaY}, \cite{M} about some attempts in this direction.

\enddocument

\end